# Point processes of exceedances for random walks in random sceneries

Ahmad Darwiche

**Abstract**

Let $\{\xi(k), k \in \mathbb{Z}\}$ be a stationary sequence of random variables and let $\{S_n, n \in \mathbb{N}_+\}$ be a transient random walk in the domain of attraction of a stable law. In the previous work [1], under conditions of type $D(u_n)$ and $D'(u_n)$ we provided a limit theorem for the maximum of the first $n$ terms of the sequence $\{\xi(S_n), n \in \mathbb{N}\}$. In this paper, under the same conditions we will see that, the limit of the process which counts the numbers of the exceedances of the form $\{\xi(S_k) > u_n\}, k \geq 1$, is a compound Poisson point process. We also deal with the so-called extremal index for the sequence $\{\xi(S_n), n \in \mathbb{N}\}$ and we discuss some weak mixing properties.

## 1 Introduction

Let $\{X_k, k \in \mathbb{N}_+\}$ be a sequence of integer-valued centered, independent identically distributed (i.i.d.) random variables and let $S_n := X_1 + \cdots + X_n$, $n \in \mathbb{N}_+$. Assume that, for any $x \in \mathbb{R}$,

$$\mathbb{P}\left(\frac{S_n}{n^{1/\alpha}} \leq x\right) \xrightarrow[n \to \infty]{} F_\alpha(x),$$

where $F_\alpha$ is the distribution function of a stable law with characteristic function given by

$$\phi(\theta) = \exp(-|\theta|^\alpha (C_1 + iC_2 \operatorname{sgn} \theta)), \quad \alpha \in (0, 2].$$

The sequence $\{S_n, n \in \mathbb{N}_+\}$ is referred to as a random walk. When $\alpha < 1$ (resp. $\alpha > 1$), it is known that the random walk is transient (resp. recurrent) [5, 6]. We denote by $R_n := \#\{S_1, \ldots, S_n\}$ the range of the random walk $\{S_n, n \in \mathbb{N}_+\}$ and by $\{Y(t), t \in \mathbb{R}\}$ the right-continuous $\alpha$-stable Lévy process with characteristic function given by $\phi(t\theta)$. In [5], it is shown that the processes $S^{(n)}(t) := \frac{S_{[nt]}}{n^{1/\alpha}}$ converges in distribution to $Y$ with respect to the Skorokhod topology ($[\cdot]$ is the integer part).

Let us recall the following result of LeGall and Rosen [6] on the asymptotic behavior of the range $R_n$.

**Theorem 1** (LeGall and Rosen). *(i) If $\alpha < 1$, then*

$$\frac{R_{[nt]}}{n} \underset{n \to \infty}{\longrightarrow} qt \quad \mathbb{P}-a.s. \tag{1}$$

*with $q := \mathbb{P}\left(S_k \neq 0, \forall k \in \mathbb{N}_+\right)$.*

*(ii) If $\alpha = 1$, then*

$$\frac{h(n) R_{[nt]}}{n} \underset{n \to \infty}{\longrightarrow} t \quad in \quad L^p(\mathbb{P}),$$

*where $h(n) := 1 + \sum_{k=1}^{n} \mathbb{P}\left(S_k = 0\right)$.*

*(iii) If $1 < \alpha \leq 2$, then for any $L \in \mathbb{N}$ and any $t_1 < \cdots < t_L$, we have*

$$\frac{1}{n^{1/\alpha}} \left(R_{[nt_1]}, \ldots, R_{[nt_L]}\right) \underset{n \to \infty}{\longrightarrow} \left(m(Y(0, t_1)), \ldots, m(Y(0, t_L))\right),$$

*in distribution, where $m$ is the Lebesgue measure on $\mathbb{R}$.*

In extreme value theory, the study of rare events is tied with the observation of abnormally high values among a series of realisations of certain variables of interest, $\{\xi(k), k \in \mathbb{Z}\}$, i.e., one is interested in events of the form $\{\xi(k) > u\}$, corresponding to the exceedance of a high threshold $u$.

Let $\{\xi(k), k \in \mathbb{Z}\}$ be a family of $\mathbb{R}$-valued i.i.d. random variables. For any sequence of real numbers $(u_n)$ and for $\tau > 0$, it is known that

$$n\mathbb{P}\left(\xi > u_n\right) \underset{n \to \infty}{\longrightarrow} \tau \iff \mathbb{P}\left(\max_{k \leq n} \xi(k) \leq u_n\right) \underset{n \to \infty}{\longrightarrow} \exp(-\tau),$$

where $\xi$ has the same distribution as $\xi(k), k \in \mathbb{Z}$ (see e.g. [4]). This property has been extended for sequences of dependent random variables satisfying the so-called conditions $D(u_n)$ and $D'(u_n)$ of Leadbetter [7, 8].

Now, assume that the threshold is of the form $u_n(x) = a_n x + b_n$ ($a_n \in \mathbb{R}$, $b_n > 0$ and $x \in \mathbb{R}$). In what follows, we consider a measure $\nu$ defined on some topological space $E$, by

$$\nu(x, \infty) := \lim_{n \to \infty} n\mathbb{P}\left(\xi > u_n(x)\right). \tag{2}$$

If $\xi$ is in the domain of attraction of the extreme value distribution $G$, then the measure $\nu$ is of the form:

$$\nu(x, \infty) = \begin{cases} x^{-\beta}, & E = (0, \infty] \quad \text{if G is a Fréchet distribution} \\ (-x)^{-\delta}, & E = (-\infty, 0] \quad \text{if G is a Weibull distribution} \\ e^{-x}, & E = (-\infty, \infty] \quad \text{if G is a Gumbel distribution} \end{cases},$$

where $\beta$ and $\delta$ are as in [4].

Our goal in this paper is to keep track of the occurrences of extreme events of the form $\{\xi(S_k) > u_n\}$ on a certain time frame and then be able to provide statements regarding its asymptotic behaviour. For this purpose, we investigate the weak convergence of the point processes of exceedances,

$$\Phi_n = \left\{ \left( \frac{\tau_k}{n}, \frac{\xi(S_{\tau_k}) - b_n}{a_n} \right), k \geq 1 \right\} \tag{3}$$

where $\tau_k = \inf\{m \in \mathbb{N}_+, \#\{S_1, \ldots, S_m\} \geq k\}$ is the instant when the random walk discovers the $k$-th site. More recently in [3], by using the Laplace functional, Franke and Saigo proved that if $\{\xi(k), k \in \mathbb{Z}\}$ is i.i.d., then the point processes of exceedances defined in (3), converges to a Poisson point process with explicit intensity measure. In this paper, we provide a new approach using the Kallenberg's theorem, to extend their results by only assuming that the sequence $\{\xi(k), k \in \mathbb{Z}\}$ satisfies conditions of type $D(u_n)$ and $D'(u_n)$. More precisely, let $\{X_n, n \in \mathbb{N}_+\}$ be as above, i.e. a sequence of centered, integer-valued i.i.d. random variables in the domain of attraction of $\alpha$-stable law. Let $\{\xi(k), k \in \mathbb{Z}\}$ be a stationary sequence of random variables independent of $\{X_k, k \in \mathbb{N}_+\}$. Stationary means that, for each $p, t \in \mathbb{N}$ and $i_1 < \ldots < i_p \in \mathbb{N}$ the joint distribution function of $\xi(i_1), \ldots, \xi(i_p)$ is equal to the joint distribution function of $\xi(i_1 + t), \ldots, \xi(i_p + t)$. In the following, the sequence $\{\xi(k), k \in \mathbb{Z}\}$ is supposed to satisfy conditions of type $D(u_n)$ and $D'(u_n)$. For the first condition, we write for each $i_1 < \ldots < i_p$ and for each $u \in \mathbb{R}$,

$$F_{i_1, \ldots, i_p}(u) = \mathbb{P}\left(\xi(i_1) \leq u, \ldots, \xi(i_p) \leq u\right).$$

**Condition $D(u_n)$** We say that $\{\xi(k), k \in \mathbb{Z}\}$ satisfies the condition $D(u_n)$ if there exist a sequence $(\alpha_{n,l})_{(n,l) \in \mathbb{N}^2}$ and a sequence $(l_n)$ of positive integers such that $\alpha_{n,l_n} \to 0$ as $n$ goes to infinity, $l_n = o(n)$, and

$$|F_{i_1, \ldots, i_p, j_1, \ldots, j_{p'}}(u_n) - F_{i_1, \ldots, i_p}(u_n) F_{j_1, \ldots, j_{p'}}(u_n)| \leq \alpha_{n,l}$$

for any integers $i_1 < \cdots < i_p < j_1 < \cdots < j_{p'}$ such that $j_1 - i_p \geq l$. Notice that the bound holds uniformly in $p$ and $p'$. Roughly, the condition $D(u_n)$ (see e.g. p29 in [10]) is a weak mixing property for the tails of the joint distributions.

The condition $D'(u_n)$ (see e.g. p29 in [10]) is a local type property and precludes the existence of clusters of exceedances. To introduce it, we consider a sequence $(k_n)$ such that

$$k_n \xrightarrow[n \to \infty]{} \infty, \quad \frac{n^2}{k_n} \alpha_{n,l_n} \xrightarrow[n \to \infty]{} 0, \quad k_n l_n = o(n), \tag{4}$$

where $(l_n)$ and $(\alpha_{n,l})_{(n,l) \in \mathbb{N}^2}$ are the same as in condition $D(u_n)$.

**Condition $D'(u_n)$** We say that $\{\xi(k), k \in \mathbb{Z}\}$ satisfies the condition $D'(u_n)$ if there exists a sequence of integers $(k_n)$ satisfying (4) such that

$$\lim_{n \to \infty} n \sum_{j=1}^{[n/k_n]} \mathbb{P}\left(\xi(0) > u_n, \xi(j) > u_n\right) = 0.$$

In the classical literature, the sequences $(\alpha_{n,l})_{(n,l)\in\mathbb{N}^2}$ and $(k_n)$ only satisfy $k_n\alpha_{n,l_n} \underset{n\to\infty}{\longrightarrow} 0$ (see e.g. (3.2.1) in [10]) whereas in (4) we have assumed that $\frac{n^2}{k_n}\alpha_{n,l_n} \underset{n\to\infty}{\longrightarrow} 0$. In this sense, the condition $\mathbf{D}'(u_n)$ as written above is slightly more restrictive than the usual condition $D'(u_n)$. Basing in the foregoing, we present our main theorem.

**Theorem 2.** *For $\alpha \leq 1$, let $u_n(x) = a_n x + b_n$ be such that $n\mathbb{P}(\xi > u_n(x)) \to \nu(x,\infty)$, as $n \to \infty$. And let $\tau_k = \inf\{m \in \mathbb{N}; \#\{S_1, \ldots, S_m\} \geq k\}$. Assume that conditions $\mathbf{D}(u_n)$ and $\mathbf{D}'(u_n)$ hold for the sequence $\{\xi(k), k \in \mathbb{Z}\}$. Then the point process*

$$N^{(n)} := \sum_k \delta_{\left(\frac{\tau_k}{n}, \frac{\xi(S_{\tau_k}) - b_{m(n)}}{a_{m(n)}}\right)}$$

*converges weakly to the Poisson point process $N$ with intensity measure $m \times \nu$, where*

$$m(n) = \begin{cases} [qn], & \text{if } \alpha < 1 \\ \left[\frac{n}{h(n)}\right], & \text{if } \alpha = 1 \end{cases}.$$

In the following, let $\mu$ be the measure defined by

$$\mu(dt, dx) = m_Y(dt) \times \nu(dx),$$

where $m_Y(t) = m(Y(0,t))$ and let $N_Y$ be the Poisson point process on $\mathbb{R}^+ \times E$ with random intensity measure $\mu$. Since $m_Y$ is a random measure, the random measure $N_Y$ is a Cox point process. The following theorem deals with the recurrent case.

**Theorem 3.** *For $\alpha > 1$, the point processes*

$$N^{(n)} := \sum_k \delta_{\left(\frac{\tau_k}{n}, \frac{\xi(S_{\tau_k}) - b_{[n^{1/\alpha}]}}{a_{[n^{1/\alpha}]}}\right)}$$

*on $\mathbb{R}^+ \times E$ converges weakly to the point process $N_Y$.*

The rest of this paper is organized as follows. In Section 2, we recall some known results, which will be used in the proofs of our main theorems. In Section 3, we calculate the extremal index for the sequence $\{\xi(S_n), n \in \mathbb{N}\}$ so-called a random walk in a random scenery. In Section 4, we discuss some weak mixing properties for the random walks in random scenery.

## 2 Proofs

In this section, we recall some known results on the convergence of the point processes and we provide some auxiliary results required for the proofs of Theorem 2 and Theorem 3, which have independent significance.

Let $E$ be a topological space. We say that $m$ is a point measure on $E$ if $m = \sum_{i=1}^{\infty} \delta_{x_i}$, where $\delta_{x_i}$ denotes the Dirac measure supported on $x_i \in E$. Consider the set $M_p(E)$ of point measures on $E$ endowed with the vague topology. A point process on $E$ is a random element of $M_p(E)$. Call a point process $\Phi$ simple if its distribution concentrates on the simple point measures of $M_p(E)$. This means

$$\mathbb{P}\left(\Phi(\{x\}) \leq 1 \quad \forall x \in E\right) = 1.$$

In order to prove our results we will use a Kallenberg's Theorem on the convergence of point processes (see Proposition 3.22 in [12]).

**Theorem 4** (Kallenberg). *Suppose $\Phi$ is a simple point process on $E$ and $\mathcal{I}$ is a basis of relatively compact open sets such that $\mathcal{I}$ is closed under finite unions and intersections and for $I \in \mathcal{I}$*

$$\mathbb{P}\left(\Phi(\partial I) = 0\right) = 1,$$

*where $\partial I$ is the boundary of $I$. Let $(\Phi_n)$ be a sequence of point processes on $E$ such that, for all $I \in \mathcal{I}$*

$$\lim_{n \to +\infty} \mathbb{P}\left(\Phi_n(I) = 0\right) = \mathbb{P}\left(\Phi(I) = 0\right)$$

*and*

$$\lim_{n \to +\infty} \mathbb{E}\left(\Phi_n(I)\right) = \mathbb{E}\left(\Phi(I)\right).$$

*Then $\Phi_n$ converges to $\Phi$ in $M_p(E)$.*

The next lemma states the independence of the sequence $\{\xi(S_n), n \in \mathbb{N}\}$ and the sequence $\{\tau_k, k \in \mathbb{N}\}$.

**Lemma 1.** *For all measurable sets $B \subset \mathbb{N}$ and $A \subset \mathbb{R}$, we have*

$$\mathbb{P}\left(\tau_k \in B, \xi(S_{\tau_k}) \in A\right) = \mathbb{P}\left(\tau_k \in B\right) \mathbb{P}\left(\xi(k) \in A\right).$$

**Proof of Lemma 1.** Since the random walk and the random scenery are independent, we have

$$\begin{aligned}
\mathbb{P}\left(\tau_k \in B, \xi(S_{\tau_k}) \in A\right) &= \sum_{m \in B} \mathbb{P}\left(\tau_k = m, \xi(S_m) \in A\right) \\
&= \sum_{m \in B} \sum_{z \in \mathbb{Z}} \mathbb{P}\left(\tau_k = m, S_m = z, \xi(z) \in A\right) \\
&= \sum_{m \in B} \sum_{z \in \mathbb{Z}} \mathbb{P}\left(\tau_k = m, S_m = z\right) \mathbb{P}\left(\xi(z) \in A\right) \\
&= \mathbb{P}\left(\tau_k \in B\right) \mathbb{P}\left(\xi(k) \in A\right).
\end{aligned}$$

$\square$

In [1], assuming that the conditions $\mathbf{D}(u_n)$ and $\mathbf{D}'(u_n)$ hold for the sequence $\{\xi(k), k \in \mathbb{Z}\}$, we have shown that the both quantities:

$$\mathbb{P}\left(\bigcap_{k \geq 1 : \frac{\tau_k}{n} \in (0,1]} \left\{\frac{\xi(S_{\tau_k}) - b_{m(n)}}{a_{m(n)}} \notin (x, \infty)\right\}\right) \quad \text{and} \quad \mathbb{E}\left(\exp\left(-\frac{R_n}{m(n)}\tau\right)\right) \qquad (5)$$

have the same behavior asymptotically, when $m(n)\mathbb{P}\left(\frac{\xi - b_{m(n)}}{a_{m(n)}} \in (x, \infty)\right)$ converges to $\tau$ as $n$ goes to infinity. This result has been proved precisely for the transient case ($\alpha < 1$) and for $m(n) = n$ (see Theorem 1 in [1]), but it remains true for the other cases of the random walk ($\alpha \geq 1$), using LeGall and Rosen's theorem according to the value of the sequence $m(n)$.

The following lemma deals with the case where the interval $(0, 1]$ (resp. $(x, \infty)$) is replaced by $(a, b]$ in (5) (resp. $A \subset \mathbb{R}$).

**Lemma 2.** *Let $\{S_n, n \in \mathbb{N}_+\}$ be as above and let $\{\xi(k), k \in \mathbb{Z}\}$ be a stationary sequence of random variables such that for any $B \subset \mathbb{R}$, $m(n)\mathbb{P}\left(\frac{\xi - a_{m(n)}}{b_{m(n)}} \in A\right) \xrightarrow[n \to \infty]{} \nu(A)$, for some sequence $(a_n) \subset \mathbb{R}$ and $(b_n) \subset \mathbb{R}_+^*$. Assume that the conditions $\mathbf{D}(u_n)$ and $\mathbf{D}'(u_n)$ hold. Then for almost all realization of $\{S_n, n \in \mathbb{N}_+\}$,*

$$\lim_{n \to \infty} \mathbb{P}\left(\bigcap_{k \geq 1 : \frac{\tau_k}{n} \in (a,b]} \left\{\frac{\xi(S_{\tau_k}) - b_{m(n)}}{a_{m(n)}} \notin A\right\}\right) - \mathbb{E}\left(\exp\left(-\frac{R_{[nb]} - R_{[na]}}{m(n)}\nu(A)\right)\right) = 0.$$

The existence of an integer $k$ such that $\tau_k$ belongs to the time interval $(a, b]$ means that the random walk $\{S_n, n \in \mathbb{N}_+\}$ visits a new site during $(a, b]$. This implies that the total number of new sites visited by $\{S_n, n \in \mathbb{N}_+\}$ during the time interval $(a, b]$ is $\sum_k \delta_{\tau_k}((a,b])$. Therefore, we have

$$\#\{k \geq 1 : \tau_k \in (a, b]\} = \sum_k \delta_{\tau_k}((a,b]) = R_{[b]} - R_{[a]}.$$

**Proof of Lemma 2.** The proof is an adaptation of the proof of Theorem 1 in [1], we replace the time interval $(0, 1]$ (resp. $(x, \infty)$) by $(a, b]$ (resp. $A \subset \mathbb{R}$) and use the fact that $\#\{k \geq 1 : \frac{\tau_k}{n} \in (a, b]\} = R_{[nb]} - R_{[na]}$.

To be more precise, let $(k_n)$, $(l_n)$ be as in (4). For $n$ large enough, let $r_n = \left[\frac{n}{k_n}\right]$. Given a realization $\{S_n, n \in \mathbb{N}_+\}$ of the random walk, we write

$$\mathcal{S}_{(na,nb]} = \{S_{\tau_k} : k \geq 1, \frac{\tau_k}{n} \in (a, b]\} \quad \text{and} \quad R_{[nb]} - R_{[na]} = \#\mathcal{S}_{(na,nb]}.$$

To capture the fact that the random scenery $\{\xi(k), k \in \mathbb{Z}\}$ satisfies the condition $\mathbf{D}(u_n)$, we construct blocks and stripes as follows. Let

$$K_n = \left[\frac{R_{[nb]} - R_{[na]}}{r_n}\right] + 1. \qquad (6)$$

There exists a unique $K_n$-tuple of subsets $B_i \subset \mathcal{S}_{(na,nb]}$, $i \leq K_n$, such that the following properties hold: $\bigcup_{j \leq K_n} B_j = \mathcal{S}_{(na,nb]}$, $\#B_i = r_n$ and $\max B_i < \min B_{i+1}$ for all $i \leq K_n - 1$. Notice that $K_n \leq \bar{k}_n$ and $\#B_{K_n} = R_{[nb]} - R_{[na]} - (K_n - 1) \cdot r_n$ almost surely (a.s.). The sets $B_j$, $j \leq K_n$, are referred to as *blocks*. For each $j \leq K_n$, we also denote by $L_j$ the family consisting of the $l_n$ largest terms of $B_j$ (e.g. if $B_j = \{x_1, \ldots, x_{r_n}\}$, with $x_1 < \cdots < x_{r_n}$, $j \leq K_n - 1$, then $L_j = \{x_{r_n - l_n + 1}, \ldots, x_{r_n}\}$). When $j = K_n$, we take the convention $L_{K_n} = \emptyset$ if $\#B_{K_n} < l_n$. The set $L_j$ is referred to as a *stripe*, and the union of the stripes is denoted $\mathcal{L}_n = \bigcup_{j \leq K_n} L_j$. Adapting Lemmas 1 and 2 of [1] in our context, we have

- $\mathbb{P}\left(\bigcap_{i \in \mathcal{S}_{(na,nb]}} \left\{\frac{\xi(i) - b_{m(n)}}{a_{m(n)}} \notin A\right\}\right) - \mathbb{P}\left(\bigcap_{i \in \mathcal{S}_{(na,nb]} \setminus \mathcal{L}_n} \left\{\frac{\xi(i) - b_{m(n)}}{a_{m(n)}} \notin A\right\}\right) \xrightarrow[n \to \infty]{} 0$,

- $\mathbb{P}\left(\bigcap_{i \in \mathcal{S}_{(na,nb]} \setminus \mathcal{L}_n} \left\{\frac{\xi(i) - b_{m(n)}}{a_{m(n)}} \notin A\right\}\right) - \prod_{j \leq K_n} \mathbb{P}\left(\bigcap_{i \in B_j \setminus \mathcal{L}_n} \left\{\frac{\xi(i) - b_{m(n)}}{a_{m(n)}} \notin A\right\}\right) \xrightarrow[n \to \infty]{} 0$,

- $\prod_{j \leq K_n} \mathbb{P}\left(\bigcap_{i \in B_j \setminus \mathcal{L}_n} \left\{\frac{\xi(i) - b_{m(n)}}{a_{m(n)}} \notin A\right\}\right) - \prod_{j \leq K_n} \mathbb{P}\left(\bigcap_{i \in B_j} \left\{\frac{\xi(i) - b_{m(n)}}{a_{m(n)}} \notin A\right\}\right) \xrightarrow[n \to \infty]{} 0$,

- $\prod_{j \leq K_n} \mathbb{P}\left(\bigcap_{i \in B_j} \left\{\frac{\xi(i) - b_{m(n)}}{a_{m(n)}} \notin A\right\}\right) - \mathbb{E}\left(\exp\left(-\frac{R_{[nb]} - R_{[na]}}{m(n)} \nu(A)\right)\right) \xrightarrow[n \to \infty]{} 0$.

The first and the third assertions mean that, asymptotically, the intersection is not affected if we remove the sites which belong to one of the stripes. Roughly, this comes from the fact that the size of the stripes is negligible compared to the size of the blocks. The second assertion is a consequence of the fact that the sequence $\{\xi(k), k \in \mathbb{Z}\}$ satisfies the condition $\mathbf{D}(u_n)$ and the last assertion is concluded by performing a calculation using conditions $\mathbf{D}(u_n)$ and $\mathbf{D}'(u_n)$. Lemma 2 follows directly from the four assertions. $\square$

According the LeGall and Rosen's Theorem and Lemma 2, we obtain the following corollary.

**Corollary 1.** *With the notation of Lemma 2, we have*

$$\lim_{n \to \infty} \mathbb{P}\left(\bigcap_{k \geq 1: \frac{\tau_k}{n} \in (a,b]} \left\{\frac{\xi(S_{\tau_k}) - b_{m(n)}}{a_{m(n)}} \notin B\right\}\right)$$

$$= \begin{cases} \exp(-(b-a)\nu(B)), & \text{if} \quad m(n) = [qn] \quad \text{and} \quad \alpha < 1 \\ \exp(-(b-a)\nu(B)), & \text{if} \quad m(n) = \frac{n}{h(n)} \quad \text{and} \quad \alpha = 1 \\ \mathbb{E}(\exp(-m_Y((a,b])\nu(B))), & \text{if} \quad m(n) = n^{\frac{1}{\alpha}} \quad \text{and} \quad \alpha > 1 \end{cases}$$

*where* $m_Y((a,b]) = m(Y(0,b)) - m(Y(0,a))$.

**Proof of Theorem 2.** We only give a detailed proof of Theorem 2 since Theorem 3 can be proved in a similar way. According to Kallenberg's Theorem, it suffices to show that for all $I \in \mathcal{I}$, we have

1. $\lim_{n \to \infty} \mathbb{E}\left(N^{(n)}(I)\right) = m \times \nu(I)$,

2. $\lim_{n \to \infty} \mathbb{P}\left(N^{(n)}(I) = 0\right) = e^{-m \times \nu(I)}$,

Without losing generality, we can assume that $I$ is of the form $I = (a, b] \times B$ with $a < b$ and $B$ is an open subset of $E$.

Firstly, we have

$$\mathbb{E}\left(N^{(n)}(I)\right) = \mathbb{E}\left(\sum_k \delta_{\left(\frac{\tau_k}{n}, \frac{\xi(S_{\tau_k}) - b_{[qn]}}{a_{[qn]}}\right)}(I)\right)$$
$$= \sum_k \mathbb{P}\left(\left(\frac{\tau_k}{n}, \frac{\xi(S_{\tau_k}) - b_{[qn]}}{a_{[qn]}}\right) \in I\right). \quad (7)$$

Now, by Lemma 1, and because the sequence $\{\xi(k), k \in \mathbb{Z}\}$ is stationary, we have

$$\mathbb{P}\left(\left(\frac{\tau_k}{n}, \frac{\xi(S_{\tau_k}) - b_{[qn]}}{a_{[qn]}}\right) \in I\right) = \mathbb{P}\left(\frac{\tau_k}{n} \in (a, b]\right) \mathbb{P}\left(\frac{\xi(S_{\tau_k}) - b_{[qn]}}{a_{[qn]}} \in B\right)$$
$$= \mathbb{P}\left(\frac{\tau_k}{n} \in (a, b]\right) \mathbb{P}\left(\frac{\xi(1) - b_{[qn]}}{a_{[qn]}} \in B\right)$$
$$= \mathbb{P}\left(\frac{\tau_k}{n} \in (a, b]\right) P_{[qn]}(B). \quad (8)$$

Combining (7) and (8), and using the fact that $\sum_k \delta_{\frac{\tau_k}{n}}((a, b]) = R_{[nb]} - R_{[na]}$, we have

$$\mathbb{E}\left(N^{(n)}(I)\right) = \sum_k \mathbb{P}\left(\frac{\tau_k}{n} \in (a, b]\right) P_{[qn]}(B)$$
$$= \mathbb{E}\left(\sum_k \delta_{\frac{\tau_k}{n}}(a, b]\right) P_{[qn]}(B)$$
$$= \mathbb{E}\left(R_{[nb]} - R_{[na]}\right) P_{[qn]}(B).$$

Moreover, by (1) and according to the dominated convergence theorem, we know that $\mathbb{E}(R_{[nb]} - R_{[na]}) \underset{n \to \infty}{\sim} nq(b - a)$. Thus we have

$$\mathbb{E}(R_{[nb]} - R_{[na]}) P_{[qn]}(B) \underset{n \to \infty}{\to} (b - a) \times \nu(B).$$

Secondly, we have

$$\mathbb{P}\left(\sum_k \delta_{\left(\frac{\tau_k}{n}, \frac{\xi(S_{\tau_k}) - b_{[qn]}}{a_{[qn]}}\right)}(I) = 0\right) = \mathbb{P}\left(\bigcap_{k \geq 1: \frac{\tau_k}{n} \in [a, b]} \left\{\frac{\xi(S_{\tau_k}) - b_{[qn]}}{a_{[qn]}} \notin B\right\}\right).$$

By Corollary 1, we obtain

$$\lim_{n \to \infty} \mathbb{P}\left(\bigcap_{k \geq 1: \frac{\tau_k}{n} \in [a, b]} \left\{\frac{\xi(S_{\tau_k}) - b_{[qn]}}{a_{[qn]}} \notin B\right\}\right) = \exp\left(-(b - a)\nu(B)\right),$$

as desired. We therefore conclude the proof of Theorem 2. □

# 3 Extremal index

The notion of extremal index appeared in the work of Loynes [9] but was formally introduced by Leadbetter in [8]. One can interpret the inverse of the extremal index as the mean number of exceedances in a cluster of exceedances.

**Definition 1.** *Let $\theta \in [0,1]$. We say that $\theta$ is the extremal index for the sequence $\{\zeta_i, i \in \mathbb{N}\}$ if, for every $\tau > 0$,*

1. *There exists a sequence $u_n \uparrow \infty$ such that $n(1 - F(u_n)) \to \tau$,*
2. *$\mathbb{P}\left(\max\{\zeta_1, \ldots, \zeta_n\} \leq u_n\right) \to \exp(-\tau\theta)$.*

Let $r_n = \left[\dfrac{n}{k_n}\right]$ where $k_n$ is as in (4). Let $M_{i,j} := \begin{cases} \max_{i \leq t \leq j} \xi(t) & i \leq j \\ -\infty & i > j \end{cases}$ and let $M_n = M_{1,n}$.
The following result is due to O'Brien and characterizes the extremal index (see Theorem 2.1. in [11]).

**Theorem 5** (O'Brien). *Let $\{\xi(k), k \in \mathbb{Z}\}$ be a stationary sequence of random variables such that $n\mathbb{P}\left(\xi > u_n\right) \xrightarrow[n \to \infty]{} \tau$, for some $\tau > 0$. Assume that the condition $D(u_n)$ holds for the sequence $\{\xi(k), k \in \mathbb{Z}\}$. Then*

$$\mathbb{P}\left(M_n \leq u_n\right) - \exp(-n\mathbb{P}\left(\xi_1 > u_n \geq M_{2,r_n}\right)) \xrightarrow[n \to \infty]{} 0.$$

In particular, if $\mathbb{P}\left(M_{2,r_n} \leq u_n | \xi_1 > u_n\right) \xrightarrow[n \to \infty]{} \theta$, then the extremal index of the stationary sequence $\{\xi(k), k \in \mathbb{Z}\}$ is equal to $\theta$. In the following theorem, we investigate the asymptotic behaviour of $\mathbb{P}\left(M_{\mathcal{S}_n} \leq u_n\right)$ under the condition $D(u_n)$, where $M_{\mathcal{S}_n} = \max_{k \in \mathcal{S}_n} \xi(k)$ with $\mathcal{S}_n = \{S_1, \ldots, S_n\}$. In the rest of the paper, we write $M_B = max_{k \in B} \xi(k)$. Although the sequence $\{\xi(S_n), n \in \mathbb{N}\}$ is not stationary, we also make explicit its extremal index. To do it, we will subdivide the set $\mathcal{S}_n$ into $K_n$ blocks in the same spirit as we did in [1]. To simplify the notation, we assume without loss of generality that the last block has the same size as the others. In other words, we suppose that $\frac{R_n}{K_n}$ is an integer.

**Theorem 6.** *Let $\{S_n, n \in \mathbb{N}_+\}$ be as above and let $\{\xi(k), k \in \mathbb{Z}\}$ be a stationary sequence of random variables such that $n\mathbb{P}\left(\xi > u_n\right) \xrightarrow[n \to \infty]{} \tau$, for some sequence $(u_n)$. Assume that the condition $D(u_n)$ holds for the sequence $\{\xi(k), k \in \mathbb{Z}\}$. Then for almost all realization of $\{S_n, n \in \mathbb{N}_+\}$,*

$$\mathbb{P}\left(M_{\mathcal{S}_n} \leq u_n\right) - \exp\left(-\sum_{j \leq K_n} \sum_{i=1}^{r_n} \mathbb{P}\left(\xi(S_{((j-1)r_n+i)}) > u_n \geq M'_{((j-1)r_n+i+1, jr_n)}\right)\right) \xrightarrow[n \to \infty]{} 0,$$

where $M'_{(i,j)} := \begin{cases} \max_{i \leq t \leq j} \xi(S_{(t)}) & i \leq j \\ -\infty & i > j \end{cases}$ and where $S_{(i)}$ is the $i$-th order statistic.

**Remark 7.**    1. *In the previous theorem we take condition $D(u_n)$ as it was written by Leadbetter.*

2. *If for almost all realization of $\{S_n, n \in \mathbb{N}_+\}$,*

$$\mu'(u_n) := \frac{1}{n} \sum_{j \leq K_n} \sum_{i=1}^{r_n} \mathbb{P}\left(M'_{((j-1)r_n+i+1,jr_n)} \leq u_n | \xi(S_{((j-1)r_n+i)}) > u_n\right) \xrightarrow[n \to \infty]{} \theta,$$

   *then the extremal index of the sequence $\{\xi(S_n), n \in \mathbb{N}\}$ is equal to $\theta$.*

3. *According to Theorem 1 in [1], if $\mathbf{D}(u_n)$ and $\mathbf{D}'(u_n)$ hold for the sequence $\{\xi(k), k \in \mathbb{Z}\}$, then for almost all realization of $\{S_n, n \in \mathbb{N}_+\}$, the limit of $\mu'(u_n)$ exists and is equal to $q$.*

**Proof.** Consider $\mathcal{S}_n = \{S_{(1)}, \ldots, S_{(R_n)}\}$ with $S_{(1)} < S_{(2)} < \cdots < S_{(R_n)}$. Let $B_j = \{S_{((j-1)r_n+1)}, \ldots, S_{(jr_n)}\}$ be a block of size $r_n$ so that $\bigcup_{j \leq K_n} B_j = \mathcal{S}_n$. Proceeding as in the proof of Theorem 1 and according to Lemma 1 in [1], we have for almost all realization of $\{S_n, n \in \mathbb{N}_+\}$,

$$\mathbb{P}\left(M_{\mathcal{S}_n} \leq u_n\right) - \prod_{j \leq K_n} \mathbb{P}\left(M_{B_j} \leq u_n\right) \xrightarrow[n \to \infty]{} 0.$$

Moreover,

$$\mathbb{P}\left(M_{B_j} \leq u_n\right) = 1 - \mathbb{P}\left(M_{B_j} > u_n\right)$$
$$= 1 - \sum_{i=1}^{r_n} \mathbb{P}\left(\xi(S_{((j-1)r_n+i)}) > u_n \geq M'_{((j-1)r_n+i+1,jr_n)}\right),$$

and

$$\mathbb{P}\left(M_{B_j} > u_n\right) \leq r_n \mathbb{P}\left(\xi > u_n\right),$$

which converges to 0 as $n$ goes to infinity. Since, for $x$ close to 0, we have $|\log(1-x)+x| \leq Cx^2$ for some $C > 0$. Thus

$$\left| \sum_{j \leq K_n} \log\left(1 - \mathbb{P}\left(M_{B_j} > u_n\right)\right) + \sum_{j \leq K_n} \mathbb{P}\left(M_{B_j} > u_n\right) \right|$$
$$\leq \sum_{j \leq K_n} \left| \log\left(1 - \mathbb{P}\left(M_{B_j} > u_n\right)\right) + \mathbb{P}\left(M_{B_j} > u_n\right) \right|$$
$$\leq C \sum_{j \leq K_n} \mathbb{P}\left(M_{B_j} > u_n\right)^2 \leq C k_n r_n^2 \mathbb{P}\left(\xi > u_n\right)^2.$$

The last term converges to 0 as $n$ goes to infinity since $k_n r_n \underset{n \to \infty}{\sim} n$, $n \mathbb{P}\left(\xi > u_n\right) \xrightarrow[n \to \infty]{} \tau$ and since $r_n \mathbb{P}\left(\xi > u_n\right) \xrightarrow[n \to \infty]{} 0$. This implies that for almost all realization of $\{S_n, n \in \mathbb{N}_+\}$,

$$\prod_{j \leq K_n} \mathbb{P}\left(M_{B_j} \leq u_n\right) - \exp\left(-\sum_{j \leq K_n} \mathbb{P}\left(M_{B_j} > u_n\right)\right) \xrightarrow[n \to \infty]{} 0.$$

This concludes the proof of Theorem 6. □

# 4 Weak mixing property for the random walk in random scenery

In this section, we give examples of sequences which satisfy (or do not satisfy) the $D(u_n)$ and $D'(u_n)$ conditions and another conditions. The following result is due to Franke and Saigo (see Propositions 2 and 3 in [3]).

**Theorem 8** (Franke and Saigo). *Let $\{\xi(k), k \in \mathbb{Z}\}$ be a family of $\mathbb{R}$−valued i.i.d. random variables which is independent of the sequence $\{X_k, k \in \mathbb{N}_+\}$. Then for $\alpha < 1$, the sequence $\{\xi(S_n), n \in \mathbb{N}\}$ satisfies the $D(u_n)$ condition and does not satisfy the $D'(u_n)$ condition with*

$$u_n := a_{[nq]}x + b_{[nq]},$$

*where $q = \mathbb{P}(S_k \neq 0, k \in \mathbb{N}_+)$.*

In the following theorem, we claim that the previous result remains true if we assume that the conditions $\mathbf{D}(u_n)$ and $\mathbf{D}'(u_n)$ hold instead of assuming that $\{\xi(k), k \in \mathbb{Z}\}$ is i.i.d..

**Theorem 9.** *Let $\{\xi(k), k \in \mathbb{Z}\}$ be a stationary sequence of random variables such that $n\mathbb{P}(\xi > u_n) \xrightarrow[n \to \infty]{} \tau$ for some sequence $(u_n)$ and $\tau > 0$. For $\alpha < 1$, assume that the conditions $\mathbf{D}(u_n)$ and $\mathbf{D}'(u_n)$ hold for the sequence $\{\xi(k), k \in \mathbb{Z}\}$. Then the sequence $\{\xi(S_n), n \in \mathbb{N}\}$ satisfies the $D(u_n)$ condition.*

**Proof.** Let $0 \leq i_1 < \cdots < i_p < j_1 < \cdots < j_{p'} \leq n$ be a family of integers, with $j_1 - i_p > l_n$ and $l_n = o(n)$. We use below the following notations,

- $R_{i_1,\ldots,i_p,j_1,\ldots,j_{p'}} = \#\{S_{i_1}, \ldots, S_{i_p}, S_{j_1}, \ldots, S_{j_{p'}}\}$,
- $R_{i_1,\ldots,i_p} = \#\{S_{i_1}, \ldots, S_{i_p}\}$,
- $R_{j_1,\ldots,j_{p'}} = \#\{S_{j_1}, \ldots, S_{j_{p'}}\}$,
- $R_{j_1,\ldots,j_{p'}}^{i_1,\ldots,i_p} = \#\{S_{i_1}, \ldots, S_{i_p}\} \cap \{S_{j_1}, \ldots, S_{j_{p'}}\} = R_{j_1,\ldots,j_{p'}} + R_{i_1,\ldots,i_p} - R_{i_1,\ldots,i_p,j_1,\ldots,j_{p'}}$.

We can assume without loss of generality that,

$$\frac{i_p}{n} \xrightarrow[n \to \infty]{} u \quad \text{and} \quad \frac{j_{p'}}{n} \xrightarrow[n \to \infty]{} v, \quad (u, v \geq 0).$$

From (1) it follows that

$$\frac{R_{1,\ldots,i_p}}{n} \xrightarrow[n \to \infty]{} qu \quad \mathbb{P}-a.s. \quad \text{and} \quad \frac{R_{i_p,i_p+1,\ldots,j_{p'}}}{n} \xrightarrow[n \to \infty]{} q(v-u) \quad \mathbb{P}-a.s..$$

Let $F'_{i_1,\ldots,i_p}(u_n) = \mathbb{P}(\xi(S_{i_1}) \leq u_n, \ldots, \xi(S_{i_p}) \leq u_n)$, then

$$|F'_{i_1,\ldots,i_p,j_1,\ldots,j_{p'}}(u_n) - F'_{i_1,\ldots,i_p}(u_n)F'_{j_1,\ldots,j_{p'}}(u_n)|$$

$$\leq \left| F'_{i_1,\ldots,i_p,j_1,\ldots,j_{p'}}(u_n) - \mathbb{E}\left(\exp\left(-\frac{R_{i_1,\ldots,i_p,j_1,\ldots,j_{p'}}}{n}\tau\right)\right) \right|$$

$$+ \left| \mathbb{E}\left(\exp\left(-\frac{R_{i_1,\ldots,i_p,j_1,\ldots,j_{p'}}}{n}\tau\right)\right) - \mathbb{E}\left(\exp\left(-\frac{R_{i_1,\ldots,i_p} + R_{j_1,\ldots,j_{p'}}}{n}\tau\right)\right) \right|$$

$$+ \left| \mathbb{E}\left(\exp\left(-\frac{R_{i_1,\ldots,i_p} + R_{j_1,\ldots,j_{p'}}}{n}\tau\right)\right) - F'_{i_1,\ldots,i_p}(u_n)F'_{j_1,\ldots,j_{p'}}(u_n) \right|.$$

Since $\mathbf{D}(u_n)$ and $\mathbf{D}'(u_n)$ conditions hold and since $n\mathbb{P}\left(\xi > u_n\right) \underset{n\to\infty}{\longrightarrow} \tau$, we obtain, by Lemma 2, that the first term of the right-hand side of the previous inequality converge to 0 as $n$ goes to infinity. To deal with the second term, we write

$$\exp\left(-\frac{R_{i_1,\ldots,i_p,j_1,\ldots,j_{p'}}}{n}\tau\right) - \exp\left(-\frac{R_{i_1,\ldots,i_p} + R_{j_1,\ldots,j_{p'}}}{n}\tau\right)$$

$$= \exp\left(-\frac{R_{i_1,\ldots,i_p} + R_{j_1,\ldots,j_{p'}}}{n}\tau\right)\left(\exp\left(\frac{R_{i_1,\ldots,i_p}^{j_1,\ldots,j_{p'}}}{n}\tau\right) - 1\right)$$

$$\leq \exp\left(\frac{R_{i_1,\ldots,i_p}^{j_1,\ldots,j_{p'}}}{n}\tau\right) - 1$$

$$\leq \exp\left(\frac{R_{1,\ldots,i_p}^{i_p,i_p+1,\ldots,j_{p'}}}{n}\tau\right) - 1.$$

Moreover,

$$\frac{R_{1,\ldots,i_p}^{i_p,i_p+1,\ldots,j_{p'}}}{n} = \frac{R_{1,\ldots,i_p}}{n} + \frac{R_{i_p,i_p+1,\ldots,j_{p'}}}{n} - \frac{R_{1,\ldots,j_{p'}}}{n}$$

converges almost surely to $qu + q(v - u) - qv = 0$ as $n \to \infty$. The third term also converges to 0 as $n$ goes to infinity, indeed,

$$\left|\mathbb{E}\left(\exp\left(-\frac{R_{i_1,\ldots,i_p} + R_{j_1,\ldots,j_{p'}}}{n}\tau\right)\right) - F'_{i_1,\ldots,i_p}(u_n)F'_{j_1,\ldots,j_{p'}}(u_n)\right|$$

$$\leq \left|\mathbb{E}\left(\exp\left(-\frac{R_{i_1,\ldots,i_p} + R_{j_1,\ldots,j_{p'}}}{n}\tau\right)\right) - \mathbb{E}\left(\exp\left(-\frac{R_{i_1,\ldots,i_p}}{n}\tau\right)\right)F'_{j_1,\ldots,j_{p'}}(u_n)\right|$$

$$+ \left|\mathbb{E}\left(\exp\left(-\frac{R_{i_1,\ldots,i_p}}{n}\tau\right)\right)F'_{j_1,\ldots,j_{p'}}(u_n) - F'_{i_1,\ldots,i_p}(u_n)F'_{j_1,\ldots,j_{p'}}(u_n)\right|$$

$$\leq \left|\mathbb{E}\left(\exp\left(-\frac{R_{i_1,\ldots,i_p}}{n}\tau\right)\left(\exp\left(-\frac{R_{j_1,\ldots,j_{p'}}}{n}\tau\right) - F'_{j_1,\ldots,j_{p'}}(u_n)\right)\right)\right|$$

$$+ \left|\mathbb{E}\left(\exp\left(-\frac{R_{i_1,\ldots,i_p}}{n}\tau\right)\right) - F'_{i_1,\ldots,i_p}(u_n)\right|$$

converges to 0 according to Lemma 2. This concludes the proof of theorem. $\square$

In [2], the authors introduce a local mixing condition which allows to express the extremal index in terms of joint distribution.

**Condition** $D^k(u_n)$. Let $\{\xi(k), k \in \mathbb{Z}\}$ be a sequence of random variables. In conjunction with the $D(u_n)$ condition, for any positive integer $k$ we say that the condition $D^k(u_n)$ holds if there exist tow sequences of integers $(k_n)$ and $(l_n)$ such that

$$k_n \to \infty, \quad k_n \alpha_{n,l_n} \to 0 \quad \text{and} \quad k_n l_n = o(n) \tag{9}$$

and
$$\lim_{n\to\infty} n\mathbb{P}\left(\xi(1) > u_n \geq M_{2,k}, M_{k+1,r_n} > u_n\right) = 0. \tag{10}$$

As mentioned in [2], (10) is implied by the condition

$$\lim_{n\to\infty} n \sum_{j=k+1}^{r_n} \mathbb{P}\left(\xi(1) > u_n \geq M_{2,k}, \xi(j) > u_n\right) = 0.$$

Observe that the last line is the condition $D'(u_n)$ if $k=1$.

Similarly to the condition $D^k(u_n)$, we introduce the $D^\infty(u_n)$ condition as follows.

**Definition 2.** *We say that the $D^\infty(u_n)$ condition holds if there exist tow sequences of integers $(k_n)$ and $(l_n)$ satisfying (9) such that*

$$\lim_{k\to\infty}\lim_{n\to\infty} n \sum_{j=k+1}^{r_n} \mathbb{P}\left(\xi(1) > u_n \geq M_{2,k}, \xi(j) > u_n\right) = 0.$$

The following result ensures that the sequence $\{\xi(S_n), n \in \mathbb{N}\}$ satisfies the above condition.

**Theorem 10.** *Let $\{\xi(k), k \in \mathbb{Z}\}$ be a stationary sequence of random variables such that $n\mathbb{P}\left(\xi > u_n\right) \underset{n\to\infty}{\longrightarrow} \tau$ for some sequence $(u_n)$ and $\tau > 0$. For $\alpha < 1$, assume that the $\mathbf{D}(u_n)$ and $\mathbf{D}'(u_n)$ conditions hold for the sequence $\{\xi(k), k \in \mathbb{Z}\}$. Then the sequence $\{\xi(S_n), n \in \mathbb{N}\}$ satisfies the $D^\infty(u_n)$ condition.*

**Proof.** Let $M'_{i,j} := \begin{cases} \max_{i\leq t \leq j} \xi(S_t) & i \leq j \\ -\infty & i > j \end{cases}$. For all integer $k$, we have

$$n \sum_{j=k+1}^{r_n} \mathbb{P}\left(\xi(S_1) > u_n \geq M'_{2,k}, \xi(S_j) > u_n\right)$$

$$= n \sum_{j=k+1}^{r_n} \mathbb{P}\left(\xi(S_1) > u_n \geq M'_{2,k}, \xi(S_j) > u_n | S_j = S_1\right) \mathbb{P}\left(S_j = S_1\right)$$

$$+ n \sum_{j=k+1}^{r_n} \mathbb{P}\left(\xi(S_1) > u_n \geq M'_{2,k}, \xi(S_j) > u_n | S_j \neq S_1\right) \mathbb{P}\left(S_j \neq S_1\right). \tag{11}$$

The first term of the right-hand side of the previous equation tends to zero when $n \to \infty$ and $k \to \infty$. Indeed

$$\mathbb{P}\left(\xi(S_1) > u_n \geq M'_{2,k}, \xi(S_j) > u_n | S_j = S_1\right) \mathbb{P}\left(S_j = S_1\right)$$
$$= \mathbb{P}\left(\xi(S_1) > u_n \geq M'_{2,k}\right) \mathbb{P}\left(S_j = S_1\right).$$

Since, $\{S_n, n \in \mathbb{N}_+\}$ is a transient random walk, this means that the sequence $\{S_n, n \in \mathbb{N}_+\}$ returns finitely in $S_1$ a.s., this gives $\sum_{j=2}^\infty \mathbb{P}\left(S_j = S_1\right) < \infty$. Therefore

$$\lim_{k\to\infty}\lim_{n\to\infty} n\mathbb{P}\left(\xi(S_1) > u_n\right) \sum_{j=k+1}^{r_n} \mathbb{P}\left(S_j = S_1\right) = 0.$$

Now, we prove that the second term of (11) goes to 0. We have

$$n \sum_{j=k+1}^{r_n} \mathbb{P}\left(\xi(S_1) > u_n \geq M'_{2,k},\ \xi(S_j) > u_n | S_j \neq S_1\right) \mathbb{P}(S_j \neq S_1)$$

$$= n \sum_{j=k+1}^{r_n} \mathbb{P}\left(\xi(S_1) > u_n \geq M'_{2,k},\ \xi(S_j) > u_n | S_j \in B^*(S_1, r_n)\right) \mathbb{P}(S_j \in B^*(S_1, r_n))$$

$$+ n \sum_{j=k+1}^{r_n} \mathbb{P}\left(\xi(S_1) > u_n \geq M'_{2,k},\ \xi(S_j) > u_n | S_j \notin B(S_1, r_n)\right) \mathbb{P}(S_j \notin B(S_1, r_n)), \quad (12)$$

where $B(S_1, r_n) := \{S \in \mathcal{S}_n;\ |S - S_1| \leq r_n\}$ and $B^*(S_1, r_n) = B(S_1, r_n) \setminus \{S_1\}$. We prove below that the last two terms of (12) converge to 0. For the first one, we write

$$n \sum_{j=k+1}^{r_n} \mathbb{P}\left(\xi(S_1) > u_n \geq M'_{2,k},\ \xi(S_j) > u_n | S_j \in B^*(S_1, r_n)\right) \mathbb{P}(S_j \in B^*(S_1, r_n))$$

$$\leq n \sum_{j=2}^{r_n} \mathbb{P}\left(\xi(0) > u_n,\ \xi(S_j - S_1) > u_n | S_j \in B^*(S_1, r_n)\right).$$

The last quantity converges to 0 because the sequence $\{\xi(k), k \in \mathbb{Z}\}$ satisfies the $\mathbf{D}'(u_n)$ condition. To deal with the second term, we write

$$n \sum_{j=k+1}^{r_n} \mathbb{P}\left(\xi(S_1) > u_n \geq M'_{2,k},\ \xi(S_j) > u_n | S_j \notin B(S_1, r_n)\right) \mathbb{P}(S_j \notin B(S_1, r_n))$$

$$\leq n \sum_{j=k+1}^{r_n} \mathbb{P}\left(\xi(S_1) > u_n,\ \xi(S_j) > u_n | S_j \notin B(S_1, r_n)\right)$$

$$\leq n \sum_{j=k+1}^{r_n} \mathbb{P}(\xi > u_n)^2 + n \sum_{j=k+1}^{r_n} |\mathbb{P}(\xi(S_1) > u_n, \xi(S_j) > u_n | S_j \notin B(S_1, r_n)) - \mathbb{P}(\xi > u_n)^2|.$$

$$(13)$$

We will also prove that the last two terms of (13) converge to 0. For the first one, we have

$$n \sum_{j=k+1}^{r_n} \mathbb{P}(\xi > u_n)^2 \leq n r_n \mathbb{P}(\xi > u_n)^2 \underset{n \to \infty}{\sim} \tau^2 \frac{r_n}{n},$$

which converges to 0 by assumption. To deal with the second term of (13), we use the $\mathbf{D}(u_n)$ condition. This gives

$$n \sum_{j=k+1}^{r_n} |\mathbb{P}(\xi(S_1) > u_n, \xi(S_j) > u_n | S_j \notin B(S_1, r_n)) - \mathbb{P}(\xi > u_n)^2|$$

$$\leq n \sum_{j=k+1}^{r_n} |\mathbb{P}(\xi(S_1) > u_n, \xi(S_j) > u_n | S_j \notin B(S_1, l_n)) - \mathbb{P}(\xi > u_n)^2|$$

$$\leq n r_n \alpha_{n, l_n}$$

$$\leq \frac{n^2}{k_n} \alpha_{n, l_n},$$

which converges to 0 as n goes to infinity according to (4). This concludes the proof of Theorem 10. □

**Corollary 2.** *Under the assumptions of Theorem 10. The sequence $\{\xi(S_n), n \in \mathbb{N}\}$ does not satisfy the $D'(u_n)$ condition.*

**Proof.** As the same spirit of (11) we have

$$n \sum_{j=2}^{r_n} \mathbb{P}\left(\xi(S_1) > u_n, \xi(S_j) > u_n\right)$$

$$= n \sum_{j=2}^{r_n} \mathbb{P}\left(\xi(S_1) > u_n, \xi(S_j) > u_n | S_j = S_1\right) \mathbb{P}\left(S_j = S_1\right)$$

$$+ n \sum_{j=2}^{r_n} \mathbb{P}\left(\xi(S_1) > u_n, \xi(S_j) > u_n | S_j \neq S_1\right) \mathbb{P}\left(S_j \neq S_1\right). \quad (14)$$

Since we have

$$1 - q \leq \sum_{j=2}^{\infty} \mathbb{P}\left(S_j = S_1\right) < \infty,$$

and since $\lim_{n \to \infty} n\mathbb{P}\left(\xi > u_n\right) = \tau$, $\tau > 0$, it follows for the first term on the right hand side of (14) that

$$\lim_{n \to \infty} n \sum_{j=2}^{r_n} \mathbb{P}\left(\xi(S_1) > u_n, \xi(S_j) > u_n | S_j = S_1\right) \mathbb{P}\left(S_j = S_1\right) \geq \tau(1 - q).$$

The last term is larger than 0 since $0 < q < 1$. □

# References

[1] N. Chenavier and A. Darwiche. Extremes for transient random walks in random sceneries under weak independence conditions. *Statistics & Probability Letters*, 158, 2020.

[2] M. R. Chernick, T. Hsing, and W. P. McCormick. Calculating the extremal index for a class of stationary sequences. *Advances in Applied Probability*, 23(4):835–850, 1991.

[3] B. Franke and T. Saigo. The extremes of random walks in random sceneries. *Advances in Applied Probability*, 41(2):452–468, 2009.

[4] B. Gnedenko. Sur la distribution limite du terme d'une série aléatoire. *Ann Math*, 44(3):423–453, 1943.

[5] H. Kesten and F. Spitzer. A limit theorem related to a new class of sel-similar processes. *Wahrscheinlichkeitsth*, 50:5–25, 1979.

[6] J.-F Le Gall and J. Rosen. The range of stable random walks. *Ann Proba*, 19(2):650–705, 1991.

[7] M. R. Leadbetter. On extreme values in stationary sequences. *Z. Wahrscheinlichkeitstheorie und Verw. Gebiete*, 28:289–303, 1973/74.

[8] M. R. Leadbetter. Extremes and local dependence in stationary sequences. *Z. Wahrsch. Verw. Gebiete*, 65(2):291–306, 1983.

[9] R. M. Loynes. Extreme values in uniformly mixing stationary stochastic processes. *Ann. Math. Statist.*, 36(3):993–999, 06 1965.

[10] V. Lucarini, D. Faranda, A.C.G.M.M. de Freitas, J.M.M. de Freitas, M. Holland, T. Kuna, M. Nicol, M. Todd, and S. Vaienti. *Extremes and Recurrence in Dynamical Systems*. Pure and Applied Mathematics: A Wiley Series of Texts, Monographs and Tracts. Wiley, 2016.

[11] G. L. O'Brien. Extreme values for stationary and Markov sequences. *Ann. Probab.*, 15(1):281–291, 1987.

[12] S. I. Resnick. *Extreme values, regular variation and point processes*, volume 4 of *Applied Probability. A Series of the Applied Probability Trust*. Springer-Verlag, New York, 1987.

Ahmad Darwiche
Laboratoire de Mathématiques Pures et Appliquées
Université du Littoral Côte d'Opale
Calais 62100
France
e-mail : ahmad.darwiche@univ-littoral.fr

Laboratoire de Mathématiques et leurs Applications de Valenciennes
Université polytechnique Hauts-de-France
59300 Valenciennes
France
e-mail : ahmad.darwiche@uphf.fr